\documentclass[11pt]{article}
\usepackage{amsmath}
\usepackage{amssymb}
\usepackage{amsthm}
\numberwithin{equation}{section}
\newtheorem{theorem}{Theorem}[section]
\newtheorem{proposition}[theorem]{Proposition}

\newtheorem{Remark}[theorem]{Remark}
\newenvironment{remark}{\begin{Remark}\rm}{\end{Remark}}
\newcommand\Proof{\noindent{\bf Proof}\quad}
\newcommand\nonu{\nonumber}
\newcommand\sa{\smallskipamount}

\newcommand\ba{\bigskipamount}
\newcommand\sLP{\\[\sa]}

\newcommand\bLP{\\[\ba]}
\newcommand\bPP{\\[\ba]\indent}
\newcommand\CC{\mathbb{C}}

\newcommand\RR{\mathbb{R}}

\newcommand\Ga{\Gamma}
\newcommand\De{\Delta}
\newcommand\al\alpha
\newcommand\ga\gamma
\newcommand\de\delta
\newcommand\si\sigma
\newcommand\LHS{left-hand side}
\newcommand\RHS{right-hand side}
\newcommand\iy\infty
\newcommand\pa\partial
\newcommand\union{\cup}
\newcommand{\hyp}[5]{\,\mbox{}_{#1}F_{#2}\!\left(
  \genfrac{}{}{0pt}{}{#3}{#4};#5\right)}
\renewcommand\Re{{\rm Re}\,}

\begin{document}
\title{On an identity by Chaundy and Bullard. I}
\author{Tom H. Koornwinder and Michael J. Schlosser}
\date{\sl Dedicated to Richard Askey on the occasion of his 75th birthday}
\maketitle
\begin{abstract}
An identity by Chaundy and Bullard writes $1/(1-x)^n$ ($n=1,2,\ldots$)
as a sum
of two truncated binomial series. This identity was rediscovered many times.
Notably, a special case was
rediscovered by I. Daubechies, while she was setting up the theory of
wavelets of compact support.
We discuss or survey many different proofs of the identity, and also its
relationship with  Gau{\ss}
hypergeometric series. We also consider the extension to complex values of
the two parameters which occur as summation bounds.
The paper concludes with a discussion of a multivariable analogue of
the identity, which was first given by Damjanovic, Klamkin and Ruehr.
We give the relationship with Lauricella hypergeometric functions
and corresponding PDE's.
The paper ends with a new proof of the multivariable case by splitting
up Dirichlet's multivariable beta integral.
\end{abstract}
\section{Introduction}
Chaundy and Bullard noted ``in passing'' the identity
\begin{equation}
\label{2}
1=(1-x)^{n+1}\sum_{k=0}^m\binom{n+k}k x^k
+x^{m+1}\sum_{k=0}^n\binom{m+k}k (1-x)^k,
\end{equation}
as a side result in their 1960 paper {\em John Smith's problem},
see \cite[p.256]{7}. Here $m,n$ are nonnegative integers.
Formula \eqref{2} can be written more succinctly as
\begin{equation}
p_{m,n}(x)+p_{n,m}(1-x)=1,
\label{7}
\end{equation}
where
\begin{equation}
p_{m,n}(x):=(1-x)^{n+1}\sum_{k=0}^m\binom{n+k}k x^k
=(1-x)^{n+1}\sum_{k=0}^m\frac{(n+1)_k}{k!}x^k
\label{8}
\end{equation}
and
\begin{equation}
(a)_k:=\begin{cases}
a(a+1)\ldots(a+k-1)&\mbox{if $k=1,2,\ldots\,$,}\\
1&\mbox{if $k=0$,}
\end{cases}
\label{119}
\end{equation}
is the {\em Pochhammer symbol}.

The Chaundy-Bullard identity \eqref{2} was rediscovered
(partially or completely) many times:
\begin{itemize}
\item
In 1971 Herrmann \cite{4} interpreted $p_{m,n}(x)$ (see \eqref{8})
as the polynomial of
degree $m+n+1$ which has a zero of order $n+1$ at $x=1$ and such that
$1-p_{m,n}(x)$ has a zero of order $m+1$ at $x=0$.
He proved this by induction
with respect to $n$ (although we think that he meant induction with respect to
$m+n$).
Essentially, although not explicitly given in \cite{4}, Herrmann's result
implies the identity \eqref{7}.
\item
The identity \eqref {2} was proposed in 1975 for the Canadian Mathematical
Olympiad (but not used there). Next it was proposed in 1976 for the
problem section of Crux Mathematicorum. 
A proof by induction by Kleiman was given there \cite{18} in 1977.
The same identity was also proposed in 1977
for the elementary problem section in the American Mathematical
Monthly by Burman. The Monthly \cite{15} gave two solutions in 1979,
one probabilistic proof by Schmitt and
one using partial fractions by Jagers.
Much later, in 1992 in the Monthly \cite{19} the probabilistic proof
was implicit in the solution of a problem about the longest expected
world series, posed in 1990 by Schuster. In 1997 the Monthly \cite{17}
had a follow-up with some non-probabilistic proofs.
\item
A two-variable analogue of the identity \eqref{2} was proposed in 1985
for the problem section in SIAM Review by Klamkin \& Ruehr. In 1986 Bosch
\& Steutel gave in this journal \cite{16} a probabilistic proof as solution.
In \cite{16} it was
also observed by Damjanovic, Klamkin \& Ruehr that there is
an $n$-variable generalization of the identity:
\begin{equation}
\sum_{i=1}^n x_i\,\sum_{k_1=0}^{a_1}\ldots\sum_{k_n=0}^{a_n}\de_{k_i,a_i}\,
\frac{(k_1+\cdots+k_n)!}{k_1!\ldots k_n!}\,x_1^{k_1}\ldots x_n^{k_n}=1
\qquad(x_1+\cdots+x_n=1).
\label{111}
\end{equation}
They gave a proof by generating functions. A probabilistic proof was
also indicated.
\item
In 1988 Daubechies \cite[Lemma 4.4]{1}, see also \cite[(6.1.7), (6.1.12)]{2},
rediscovered the case $m=n$ of \eqref{2}.
This identity was a crucial step
for her in order to arrive at the form of the function $m_0(\xi)$
which is associated with the wavelets of compact support named after her.
Her proof in \cite{2} was essentially the same as Jagers' proof in \cite{15},
but she referred to B\'ezout's identitity. Next Zeilberger \cite{13} in 1993
gave a probabilistic proof of Daubechies' case $m=n$ of \eqref{2} and,
unaware of \cite{16}, he stated the case $a_1=\ldots=a_n$ of \eqref{111}
and indicated a probabilistic proof.
\item
Multiplication of both sides of \eqref{2} by
$(1-x)^{-n-1}$ gives
\begin{equation}
\label{1}
(1-x)^{-n-1}=\sum_{k=0}^m\frac{(n+1)_k}{k!}x^k
+x^{m+1}\sum_{k=0}^n\frac{(m+1)_k}{k!}(1-x)^{k-n-1}.
\end{equation}
This identity is the case $m=0$ of the identity at the end of section 8
in Vid\=unas \cite{8}, where it is given as a three-term identity
for three Gau{\ss} hypergeometric functions satisfying the same
hypergeometric differential equation in the most degenerate case
(trivial monodromy group).
\item
It was also essentially identity \eqref{1} which was rediscovered by
Pieter de Jong (Netherlands), who is studying the mathematical foundations of
architecture. It was by his communication to the first author in 2007
that we first became aware of this identity.
See also de Jong's manuscript \cite{5}.
\end{itemize}

It is not without precedents in mathematics, in particular in special function
theory, that a relatively elementary result is rediscovered and published
many times. We think that for \eqref{2} the elegance and unexpectedness
of the identity arose people's interest again and again, as it did with us.
Why then spend again a publication on it? First it seems useful to survey
all earlier (as far as we know now) occurrences and approaches. Second,
we can offer some approaches which did not yet occur, notably the
approach by splitting up the beta integral and the context of the
Gau{\ss} hypergeometric function. Third the possible generalizations are
interesting.
As we mentioned, the $n$-variable generalizations already occurred, but
we can offer yet unexplored aspects of it.
There are also various analogues and
generalizations of \eqref{2} in the $q$-case, which we will present in a
forthcoming paper.

The following sections present or survey many different proofs of \eqref{2}.
The proof in section~\ref{33} is the original proof by Chaundy and Bullard
\cite{7}, and its slight variations by Daubechies \cite[\S6.1]{2}
and Jagers \cite{15} are also discussed there.
The proof in section~\ref{16} is by induction with respect to $m+n$.
The proof in section \ref{15} is by repeated differentiation of \eqref{12}
(a method suggested by Pieter de Jong in an earlier version of \cite{5}). 
The proof in section \ref{42} uses generating functions. It was communicated
to us by Helmut Prodinger and it is also a one-variable specialization
of a proof in \cite{16}.
Another proof of combinatorial flavour in section \ref{107} uses weighted
lattice paths, and becomes by specialization a probabilistic proof
(different in formulation but in essence the same as many earlier proofs
which appeared).
Section \ref{17} gives a proof by splitting up
the beta integral.
In section~\ref{19} we consider the extension of \eqref{2} to complex values
of $m,n$.
In section~\ref{43}
we observe that
the three terms in the identity \eqref{1} all solve a very special (degenerate)
case of the hypergeometric differential equation. Next
we obtain \eqref{1} as a limit case of a more general three-term identity
for hypergeometric functions, where the three terms all solve a hypergeometric
differential equation.
\par
Section \ref{108} starts our discussion of the
multivariable analogue \eqref{111}, which was first obtained in \cite{16}.
A connection with Appell and Lauricella hypergeometric functions is made.
In section \ref{109} a partial differential
equation satisfied by all terms in this multivariable analogue is given.
Finally we give in section~\ref{110}
a proof of \eqref{111} by splitting up
Dirichlet's multivariable beta integral
(generalizing the approach in section \ref{17}).
\bLP
{\bf Acknowledgements}\quad
We thank Pieter de Jong for showing us his identity and suggesting the proof
in section \ref{15}, Peter Paule
for suggesting the approach with a multivariable beta integral in section
\ref{110}, and Mike Keane, Oleg Ogievetsky and Helmut Prodinger for
helpful remarks.

The second author was partly supported by FWF Austrian Science Fund
grants \hbox{P17563-N13}, and S9607 (the latter is part
of the Austrian National Research Network
``Analytic Combinatorics and Probabilistic Number Theory'').
\bLP
{\bf Notation}\quad
The {\em Gau{\ss} hypergeometric series} (see \cite[Ch.~2]{3})
is defined by
\begin{equation}
\hyp21{a,b}{c}{z}:=\sum_{k=0}^\iy\frac{(a)_k\,(b)_k}{(c)_k\,k!}\,z^k\qquad
(z,a,b,c\in\CC,\;|z|<1,\;c\notin\{0,-1,-2,\ldots\}),
\label{116}
\end{equation}
where the Pochhammer symbol is given by \eqref{119}.
In the {\em terminating} case we have
\begin{equation}
\hyp21{-n,b}{c}{z}:=\sum_{k=0}^n\frac{(-n)_k\,(b)_k}{(c)_k\,k!}\,z^k
\quad(z,b,c\in\CC,\;n=0,1,2,\ldots,\;c\ne0,-1,\ldots,-n+1).
\label{117}
\end{equation}
\section{Chaundy \& Bullard's original proof}
\label{33}
Fix $m$ and $n$. By the binomial theorem we have
\begin{equation}
(x+y)^{m+n+1}=y^{n+1}P_{m,n}(x,y)+x^{m+1}P_{n,m}(y,x),
\label{31}
\end{equation}
where
\begin{equation}
P_{m,n}(x,y):=\sum_{k=0}^m\binom{m+n+1}k x^k y^{m-k}.
\label{112}
\end{equation}
is a homogeneous polynomial of degree $m$. Put $y:=1-x$. Then
\begin{equation}
1=(1-x)^{n+1}P_{m,n}(x,1-x)+x^{m+1}P_{n,m}(1-x,x),
\label{28}
\end{equation}
and multiplication by $(1-x)^{-n-1}$ yields
\begin{equation}
(1-x)^{-n-1}=P_{m,n}(x,1-x)+x^{m+1}(1-x)^{-n-1}P_{n,m}(1-x,x).
\label{27}
\end{equation}
Expand both sides of \eqref{27} as a power series in $x$,
convergent for $|x|<1$. Then $P_{m,n}(x,1-x)$ is a polynomial of degree $\le m$
in $x$ and all terms in
the power series of $x^{m+1}(1-x)^{-n-1}P_{n,m}(1-x,x)$ have degree $\ge m+1$.
Hence $P_{m,n}(x,1-x)$ equals the power
series of $(1-x)^{-n-1}$ truncated after
the term of $x^m$, i.e.,
\begin{equation}
P_{m,n}(x,1-x)=\sum_{k=0}^m\frac{(n+1)_k}{k!}\,x^k.
\label{29}
\end{equation}
Then substitution of \eqref{29}
in \eqref{28} proves
\eqref{2}, and its homogeneous form
\begin{equation}
(x+y)^{m+n+1}
=y^{n+1}\sum_{k=0}^m\frac{(n+1)_k}{k!} x^k(x+y)^{m-k}
+x^{m+1}\sum_{k=0}^n\frac{(m+1)_k}{k!} y^k(x+y)^{n-k}.
\label{69}
\end{equation}

Note that, conversely, \eqref{2} implies \eqref{29}, i.e., the equality
\begin{equation}
\sum_{k=0}^m\binom{m+n+1}k x^k (1-x)^{m-k}=
\sum_{k=0}^m\binom{n+k}k\,x^k.
\label{114}
\end{equation}
Indeed, compare \eqref{1} with \eqref{27}. It was essentially
this identity \eqref{114} which was also stated by Guenther \cite{14}
in a comment to the solution of a problem in the Monthly. He pointed out
many relationships of this identity with the binomial and negative
binomial distribution, including a probabilistic proof.
\begin{remark}
\label{56}
We can rewrite \eqref{114} as
\begin{equation*}
\sum_{k=0}^m\frac{(n+1)_k}{k!}\,x^k=
(1-x)^m\sum_{k=0}^m\frac{(-m-n-1)_k}{k!}\,\Bigl(\frac x{x-1}\Bigr)^k.
\end{equation*}
In terms of terminating Gau{\ss} hypergeometric series \eqref{117}
this can be written as
\begin{equation}
\hyp21{-m,n+1}{-m}x=(1-x)^m\,\hyp21{-m,-m-n-1}{-m}{\frac x{x-1}},
\end{equation}
which is the limit case $a:=-m$, $b:=n+1$, $c\to-m$
of Pfaff's transformation formula
\begin{equation}
\hyp21{a,b}cx=(1-x)^{-a}\,\hyp21{a,c-b}c{\frac x{x-1}},
\label{48}
\end{equation}
see \cite[(2.2.6)]{3}.
\end{remark}
\begin{remark}
For the proof of \eqref{2} from \eqref{28} as given above, one may start
with a weaker form of \eqref{28} as follows:
\begin{equation}
1=(1-x)^{n+1} q_{m,n}(x)+x^{m+1} r_{m,n}(x),
\label{115}
\end{equation}
where $q_{m,n}(x)$ and $r_{m,n}(x)$ are polynomials of degree $\le m$
respectively
$\le n$, so not yet necessarily explicitly given. Since $(1-x)^{n+1}$ and
$x^{m+1}$ are polynomials without common zeros of degree $n+1$ respectively
$m+1$, we can recognize \eqref{115} as a {\em B\'ezout identity},
where $q_{m,n}(x)$ and $r_{m,n}(x)$ will uniquely exist as polynomials
of precise degree $m$ respectively $n$ (see for instance
\cite[Theorem 6.1.1]{2}). It was in this way that
Daubechies \cite[\S6.1]{2}, in the symmetric case $m=n$, proved \eqref{2}.
Also note that the symmetry of
\eqref{115} together with the uniqueness and degree properties
of $q_{m,n}(x)$ and $r_{m,n}(x)$ already imply that
$q_{m,n}(x)=r_{n,m}(1-x)$, without explicit computation.

Equivalent to the B\'ezout identity approach, \eqref{115} can be seen as
a partial fraction decomposition
\[
\frac1{x^{m+1}(1-x)^{n+1}}=\frac{q_{m,n}(x)}{x^{m+1}}+
\frac{r_{m,n}(x)}{(1-x)^{n+1}}
\]
with $q_{m,n}(x)$ and $r_{m,n}(x)$ of degree $\le m$
respectively $\le n$,
cf.\ Jagers' proof of \eqref{2} in \cite{15}.
\end{remark}
\section{A proof by induction}
\label{16}
The following proof by induction was essentially given earlier by
Kleiman \cite{18} and, for $m=n$,
by Daubechies \cite[Lemma 4.4]{1}.
We have to prove \eqref{7}, with $p_{m,n}(x)$ given by \eqref{8}.
First note that \eqref{7} holds for $n=0$, and hence, by symmetry, also
for $m=0$. Indeed,
\begin{equation*}
p_{m,0}(x)=(1-x)\sum_{k=0}^m x^k=(1-x)\,\frac{1-x^{m+1}}{1-x}=1-x^{m+1},
\end{equation*}
and $p_{0,m}(1-x)=x^{m+1}$, so $p_{m,0}(x)+p_{0,m}(1-x)=1$.

Now we prove \eqref{7} by induction with respect to $m+n$. We already
saw that it holds for $m+n=0$, i.e., for $m=n=0$. Now suppose that
\eqref{7} holds for all $m,n$ with $m+n=N-1$. Let $m+n=N$.
Then we already proved \eqref{7} if $m=0$ or $n=0$, so we may assume that
$m,n>0$. Substitute the recurrence relation
\begin{equation*}
\binom{n+k}k=\binom{n+k-1}{k-1}+\binom{n+k-1}k
\end{equation*}
for binomial coefficients into \eqref{8}. Then
\begin{align}
p_{m,n}(x)&=
(1-x)^{n+1}\sum_{k=1}^m\binom{n+k-1}{k-1}x^k
+(1-x)^{n+1}\sum_{k=0}^m\binom{n+k-1}k x^k
\nonu\\
\noalign{\allowbreak}
&=x(1-x)^{n+1}\sum_{l=0}^{m-1}\binom{n+l}l x^l+(1-x) p_{m,n-1}(x)
\nonu\\
\noalign{\allowbreak}
&=x p_{m-1,n}(x)+(1-x) p_{m,n-1}(x).
\label{9}
\end{align}
Hence
\begin{equation}
p_{n,m}(1-x)=x p_{n,m-1}(1-x)+(1-x) p_{n-1,m}(1-x).
\label{10}
\end{equation}
Adding \eqref{9} and \eqref{10} gives
\begin{align*}
p_{m,n}(x)+p_{n,m}(1-x)&=x\bigl(p_{m-1,n}(x)+p_{n,m-1}(1-x)\bigr)
+(1-x)\bigl(p_{m,n-1}(x)+p_{n-1,m}(1-x)\bigr)\\
&=x+(1-x)=1
\end{align*}
by induction. This completes the proof of \eqref{7}.
\section{A proof by repeated differentiation}
\label{15}
Here we give a proof of \eqref{1} which was sketched by Pieter de Jong in
an earlier version of \cite{5}.
First note that the case $n=0$ of \eqref{1} is evident. It is essentially
the summation formula for the terminating geometric series:
\begin{equation}
\frac1{1-x}=1+x+x^2+\cdots+x^m+\frac{x^{m+1}}{1-x}\,.
\label{12}
\end{equation}
Now we can prove \eqref{1} by induction with respect
to $n$ (for each $n$ for general $m$).
For $n=0$ we have \eqref{12}, which is evident.
Apply the operator $(n+1)^{-1}\,d/dx$ to both sides of \eqref{1}.
The \LHS\ becomes $(1-x)^{-n-2}$, the \RHS\ becomes
\begin{equation*}
\sum_{k=0}^{m-1}\frac{(n+2)_k}{k!}\,x^k
+x^m\sum_{k=0}^n\frac{m+1}{n+1}\,\frac{(m+1)_k}{k!}(1-x)^{k-n-1}
+x^{m+1}\sum_{k=0}^n\frac{n-k+1}{n+1}\,\frac{(m+1)_k}{k!}(1-x)^{k-n-2}\,.
\end{equation*}
So we will have proved \eqref{1} with $n$ replaced by $n+1$ and $m$ replaced
by $m-1$ if we can show that
\begin{equation}
x^m\,\sum_{k=0}^n\frac{m+1}{n+1}\,\frac{(m+1)_k}{k!}(1-x)^{k-n-1}
+x^{m+1}\,\sum_{k=0}^n\frac{n-k+1}{n+1}\,\frac{(m+1)_k}{k!}(1-x)^{k-n-2}
\label{13}
\end{equation}
is equal to
\begin{equation}
x^m\,\sum_{k=0}^{n+1}\frac{(m)_k}{k!}(1-x)^{k-n-2}\,.
\label{14}
\end{equation}
In order to show this, rewrite $x^{m+1}$ in the second term of
\eqref{13} as $x^{m+1}=x^m-x^m(1-x)$, by which \eqref{13} becomes
\begin{align*}
&x^m\sum_{k=0}^n\frac{m-n+k}{n+1}\,\frac{(m+1)_k}{k!}(1-x)^{k-n-1}
+x^m\,\sum_{k=0}^n\frac{n-k+1}{n+1}\,\frac{(m+1)_k}{k!}(1-x)^{k-n-2}\\
\noalign{\allowbreak}
&=
x^m\sum_{k=1}^{n+1}\frac{m-n+k-1}{n+1}\,\frac{(m+1)_{k-1}}{(k-1)!}(1-x)^{k-n-2}
+x^m\,\sum_{k=0}^n\frac{n-k+1}{n+1}\,\frac{(m+1)_k}{k!}(1-x)^{k-n-2}\\
\noalign{\allowbreak}
&=x^m(1-x)^{-n-2}+x^m\,\sum_{k=1}^n\frac{(m)_k}{k!}(1-x)^{k-n-2}+
\frac{(m)_{n+1}}{(n+1)!}x^m(1-x)^{-1},
\end{align*}
which equals \eqref{14}. This completes the induction step.
\section{A proof by generating functions}
\label{42}
In \cite{16} a proof by generating functions was given for
the $n$-variable generalization \eqref{111} of \eqref{2}.
Of course this proof can be specialized to a proof by generating functions
of \eqref{2}.
Such a proof of \eqref{2} was also communicated to us by Helmut Prodinger,
independently from \cite{16}. Because the one-variable case is more simple,
we give the proof here.

Fix $x\in(0,1)$. For $u,v\in(0,1)$ let
\begin{equation}
f(u,v;x):=\sum_{m,n\ge0}u^mv^n(1-x)^{n+1}\sum_{k=0}^m\binom{n+k}k x^k.
\label{44}
\end{equation}
Then
\begin{equation}
f(v,u;1-x)=\sum_{m,n\ge0}u^mv^nx^{m+1}\sum_{k=0}^n\binom{m+k}k (1-x)^k.
\label{45}
\end{equation}
{}From \eqref{44} we have
\begin{align*}
f(u,v;x)&=
\frac1{1-u}\,\sum_{n\ge0}v^n(1-x)^{n+1}\sum_{k\ge0}\binom{n+k}k(ux)^k\\
&\;=\frac1{1-u}\,\sum_{n\ge0}v^n(1-x)^{n+1}\,\frac1{(1-ux)^{n+1}}\\
&\;=\frac{1-x}{(1-u)(1-ux)}\,\frac1{1-\frac{v(1-x)}{1-ux}}\\
&\;=\frac{1-x}{1-u}\,\frac1{1-ux-v(1-x)}\,.
\end{align*}
Hence
\[
f(v,u;1-x)=\frac x{1-v}\,\frac1{1-ux-v(1-x)}
\]
and
\[
f(u,v;x)+f(v,u;1-x)=\frac1{1-ux-v(1-x)}\,\Bigl(\frac{1-x}{1-u}+\frac x{1-v}
\Bigr)=\frac1{(1-u)(1-v)}.
\]
So
\[
f(u,v;x)+f(v,u;1-x)=\sum_{m,n\ge0}u^mv^n,
\]
and combined with \eqref{44}, \eqref{45} this yields \eqref{2}
by taking the coefficient of $u^mv^n$.
\section{A proof by weighted lattice paths}
\label{107}
Consider all lattice paths from $(0,0)$ to $(m+1,n+1)$
in the planar integer lattice (using only unit
east and north steps).
Such a path $P$ consists of $m+n+2$ successive unit steps $s_k(P)$
($k=1,2,\ldots,m+n+2$). Let $P_k$ be the path $P$ terminated after $k$ steps.
The weight $w(P)$ of a path $P$ is defined to be the
product of the weight of the respective steps $s$ of
the path, i.e., $w(P)=\prod_{s\in P}w(s)$.
Define the weight function $w$ as follows:
\begin{align*}
w\bigl((i,j)\to(i+1,j)\bigr)&:=
\begin{cases}
x&(j<n+1),\\
x+y&(j=n+1),
\end{cases}\\
w\bigl((i,j)\to(i,j+1)\bigr)&:=
\begin{cases}
y&(i<m+1),\\
x+y&(i=m+1).
\end{cases}
\end{align*}
Since for each $k\in\{1,\ldots,m+n+2\}$ and for each truncated path $P_{k-1}$
we have
\[
\sum_{P_k;\;\;P_k\backslash s_k(P_k)=P_{k-1}}w(s_k(P_k))=x+y,
\]
we find by induction that $\sum_{P_k} w(P_k)=(x+y)^k$, and hence
\begin{equation}
\sum_P w(P)=(x+y)^{m+n+2}.
\label{70}
\end{equation}

On the other hand each path $P$ ends either with a vertical step or with a
horizontal step.
Consider first the paths which end with a vertical step.
Then the last horizontal step will be $(m,k)\to(m+1,k)$ for some
$k\in\{0,1,\ldots,n\}$. For given $k$ all such paths have weight
$x^{m+1}y^k(x+y)^{n-k+1}$ and the number of such paths is
$\binom{m+k}k$. Hence the sum of the weights of all paths which end with a
vertical step equals
\begin{equation}
x^{m+1}\sum_{k=0}^n\binom{m+k}k y^k(x+y)^{n-k+1}.
\label{71}
\end{equation}
Similarly, the sum of the weights of all paths which end with a
horizontal step equals
\begin{equation}
y^{n+1}\sum_{k=0}^m\binom{n+k}k x^k(x+y)^{m-k+1}.
\label{72}
\end{equation}
Since $\eqref{70}=\eqref{71}+\eqref{72}$, we have obtained \eqref{69} with
both sides multiplied by $x+y$.
\begin{remark}
\label{118}
For $0\le x\le 1$ and $y:=1-x$ we can give a probabilistic interpretation
of the
above results. Now consider all lattice paths from $(0,0)$ to $(m+1,n+1)$,
where each following step has probability 1 if there is only one possible
step, and otherwise probability $x$ if the step is horizontal and $1-x$
if it is vertical. Then \eqref{71} gives the probability that the last
step is vertical and \eqref{72} the probability that the last
step is horizontal, and the sum of both probabilities is necessarily 1.
Thus we have a probabilistic proof of \eqref{2}.

Many probabilistic proofs of \eqref{2} were earlier given, see \cite{15},
\cite{19}, \cite{13}. They are essentially all equivalent to the one
given in the previous paragraph. Zeilberger's \cite{13} proof (phrased by him
for $m=n$) is particularly succinct. For general $m,n$ it reads as follows:
\sLP
Toss a coin (with $Pr(head)=x$) until reaching $m+1$ heads or $n+1$ tails.
Then equate the probability, 1, of finishing after at most $m+n+1$ tossings
with the sum of the probabilities of all the final outcomes. This yields
\eqref{2}.
\end{remark}
\section{A proof using the beta integral}
\label{17}
By the evaluation of the beta integral
(see \cite[\S1.1]{3})
we have for $x\in(0,1)$:
\begin{align}
1&=\frac{(m+n+1)!}{m!\,n!} \int_0^1 t^m(1-t)^n\,dt
\nonu\\
&=\frac{(m+n+1)!}{m!\,n!} \int_0^x t^m(1-t)^n\,dt+
\frac{(m+n+1)!}{m!\,n!} \int_x^1 t^m(1-t)^n\,dt
\nonu\\
&=\frac{(m+n+1)!}{m!\,n!} \int_0^x t^m(1-t)^n\,dt+
\frac{(m+n+1)!}{m!\,n!} \int_0^{1-x} t^n(1-t)^m\,dt.
\label{3}
\end{align}
Then \eqref{2} will follow from \eqref{3}
if we can prove that
\begin{equation}
\frac{(m+n+1)!}{m!\,n!} \int_0^x t^m(1-t)^n\,dt=
x^{m+1}\sum_{k=0}^n\frac{(m+1)_k}{k!}(1-x)^k.
\label{4}
\end{equation}
But \eqref{4} follows by the string of equalities
\begin{multline*}
\int_0^x t^m(1-t)^n\,dt
=x^{m+1}\int_0^1s^m\bigl(1-s+s(1-x)\bigr)^n\,ds\\
=x^{m+1}\,\sum_{k=0}^n\binom nk (1-x)^k \int_0^1 s^{m+k}(1-s)^{n-k}\,ds
=\frac{m!\,n!\,x^{m+1}}{(m+n+1)!}\,\sum_{k=0}^n\frac{(m+1)_k}{k!}(1-x)^k.
\end{multline*}

The integral on the \LHS\ of \eqref{4} is an {\em incomplete beta function},
which is usually expressed as a hypergeometric function \eqref{116}
(see \cite[\S2.5.3, p.87]{6}, also for $m$, $n$ complex with $\Re m>-1$):
\begin{equation}
B_x(m+1,n+1):=
\int_0^x t^m(1-t)^n\,dt
=\frac1{m+1}\,x^{m+1}\,\hyp21{-n,m+1}{m+2}x
\quad(x\in(0,1)).
\label{35}
\end{equation}
The proof of \eqref{35} is
by binomial expansion of $(1-t)^n$. Then \eqref{3} takes the form
\begin{multline}
1=\frac{\Gamma(n+m+2)}{\Gamma(m+1)\Gamma(n+1)}\,B_x(m+1,n+1)+
\frac{\Gamma(n+m+2)}{\Gamma(m+1)\Gamma(n+1)}\,B_{1-x}(n+1,m+1)\\
(x\in(0,1),\;m,n\in\CC,\;\Re m, \Re n>-1).
\label{122}
\end{multline}

The \RHS\ of \eqref{4} can be written as
a terminating hypergeometric series \eqref{117}.
Then combination of \eqref{4} and \eqref{35} yields
\begin{equation}
x^{m+1}\,\hyp21{-n,m+1}{-n}{1-x}
=\frac{(m+n+1)!}{(m+1)!\,n!}\,x^{m+1}\,\hyp21{-n,m+1}{m+2}x.
\label{23}
\end{equation}
Alternatively, \eqref{23} can be proved as
the limit case for $c\to -n$ of {\em Pfaff's identity}
\begin{equation}
\hyp21{-n,b}c{1-x}=\frac{(c-b)_n}{(c)_n}\,\hyp21{-n,b}{b-c-n+1}x\qquad
\mbox{($n$ nonnegative integer)},
\label{40}
\end{equation}
see \cite[(2.3.14)]{3}.
\section{Extension of the identity to non-integer $m$ and $n$}
\label{19}
By \eqref{4}, \eqref{35} the formula
\begin{equation}
p_{m,n}(x)=\frac{\Ga(m+n+2)}{\Ga(m+1)\,\Ga(n+1)}\,B_{1-x}(n+1,m+1)\quad
(x\in(0,1),\;m,n\in\CC,\;\Re n>-1)
\label{120}
\end{equation}
extends \eqref{8} to non-integer values of $m,n$. Then, by \eqref{122},
the identity \eqref{7} holds for $x\in(0,1)$ and $m,n\in\CC$ with
$\Re m, \Re n>-1$ if $p_{m,n}(x)$ is given by \eqref{120}.
Moreover, by Carlson's theorem (see for instance
Titchmarsh \cite[\S5.81]{20}) this is the unique extension
\begin{equation*}
p_{m,n}(x)+q_{m,n}(x)=1\qquad
(x\in(0,1),\;m,n\in\CC,\;\Re m,\Re n>-1)
\end{equation*}
of \eqref{7} such that for nonnegative integer $m,n$ we have
$p_{m,n}(x)=q_{n,m}(1-x)$ given by \eqref{4}, $p_{m,n}(x)$ and $q_{m,n}(x)$
are analytic in $m,n$ for $\Re m,\Re n>-1$ with $x$ fixed, and
$p_{m,n}(x)$ and $q_{m,n}(x)$ satisfy, for some $c\in(0,\pi)$,
estimates $O(e^{c|m|})$ and $O(e^{c|n|})$ as $\Re m,\Re n\ge0$.
Indeed, fix $m$ with $\Re m>-1$ and $x\in(0,1)$. Then, in the \RHS\ of
\eqref{120} we have for $\Re n\ge0$:
\[
|B_{1-x}(n+1,m+1)|\le B_{1-x}(\Re n+1,\Re m+1)\le
\int_0^{1-x} (1-t)^{\Re m}\,dt
=\frac{(1-x)^{\Re m+1}}{\Re m+1}
\]
and (as a consequence of the asymptotic formula for $\Gamma(z)$,
see \cite[Theorem 1.4.1]{3})
\[
\left|\frac{\Ga(m+n+2)}{\Ga(m+1)\,\Ga(n+1)}\right|=O(|n|^{\Re m+1}).
\]
Hence, for $x,m$ fixed as above, the \RHS\ of \eqref{120} is
$O(e^{c|n|})$ as $\Re n\ge0$ for arbitrary small $c>0$.
We can estimate the other cases in a similar way.

Alternatively, we may write \eqref{8} as
\begin{equation}
p_{m,n}(x)=(1-x)^{n+1}\sum_{k=0}^m \frac{\Ga(n+k+1)}{\Ga(n+1)\,\Ga(k+1)}\,x^k.
\label{121}
\end{equation}
Consider \eqref{121} for $x\in(0,1)$ and $n\in\CC$ with $\Re n>-1$,
and then try on it the
fractional extension of finite sums proposed
by M\"uller \& Schleicher \cite{21}, \cite{22}. Since for $k\in\CC$
with $\Re k\ge0$ we have
\begin{equation*}
f(k):=\frac{\Ga(n+k+1)}{\Ga(n+1)\,\Ga(k+1)}\,x^k(1-x)^{n+1}=o(1)\quad{\rm as}\;
\Re k\to\iy,
\end{equation*}
their recipe of fractional extension
(see \cite[(10)]{21}, \cite[top of p.5]{22}) of the sum $\sum_{k=0}^m f(k)$
is
\begin{align}
p_{m,n}(x)&=\sum_{k=0}^\iy(f(k)-f(k+m+1))
\nonu\\
&=(1-x)^{n+1}\sum_{k=0}^\iy\frac{(n+1)_k}{k!}\,x^k-(1-x)^{n+1}
\sum_{k=0}^\iy \frac{\Ga(n+m+k+2)}{\Ga(n+1)\,\Ga(m+k+2)}\,x^{k+m+1}
\nonu\\
&=1-\frac{\Gamma(n+m+2)}{\Gamma(m+2)\Gamma(n+1)}\,x^{m+1}(1-x)^{n+1}\,
\hyp21{n+m+2,1}{m+2}x
\nonu\\
&=1-\frac{\Gamma(n+m+2)}{\Gamma(m+2)\Gamma(n+1)}\,x^{m+1}\,
\hyp21{-n,m+1}{m+2}x
\nonu\\
&=1-\frac{\Gamma(n+m+2)}{\Gamma(m+1)\Gamma(n+1)}\,B_x(m+1,n+1).
\label{123}
\end{align}
Here we used Euler's transformation formula \cite[(2.2.7)]{3} and
\eqref{35}.
Note that for $x\in(0,1)$ and $m,n\in\CC$ with $\Re m,\Re n>-1$ the extension
of $p_{m,n}(x)$ defined by \eqref{120} is equal to the extension defined
by \eqref{123}. This equality is given by \eqref{122}. Curiously, this equality
is also the extension of the identity \eqref{7}.
\section{Three-term hypergeometric identities}
\label{43}
We can write \eqref{1} as
\begin{equation}
u_3=u_1+u_2,
\label{37}
\end{equation}
where
\begin{align}
u_1(x)&:=\sum_{k=0}^m\frac{(n+1)_k}{k!}x^k
=\hyp21{-m,n+1}{-m}x,
\label{46}\\
u_2(x)&:=x^{m+1}\sum_{k=0}^n\frac{(m+1)_k}{k!}(1-x)^{k-n-1}
=x^{m+1}(1-x)^{-n-1}\,\hyp21{-n,m+1}{-n}{1-x}
\nonu\\
&\,=\frac{(m+n+1)!}{(m+1)!\,n!}\,x^{m+1}(1-x)^{-n-1}\,\hyp21{m+1,-n}{m+2}x,
\label{25}\\
u_3(x)&:=(1-x)^{-n-1}=(1-x)^{-n-1}\,\hyp21{0,-m-n-1}{-n}{1-x}.
\label{47}
\end{align}
Here the third identity in \eqref{25} is \eqref{23}.
Now $u_1$, $u_2$ and $u_3$ are three different solutions of
the hypergeometric differential equation
\begin{equation}
x(1-x)u''(x)-\bigl((n+2)x+m(1-x)\bigr)u'(x)+m(n+1)u(x)=0.
\label{39}
\end{equation}
Indeed, consider Kummer's 24 solutions of the hypergeometric equation
\begin{equation}
z(1-z)u''(z)+\bigl(c-(a+b+1)z\bigr)u'(z)-abu(z)=0
\label{26}
\end{equation}
in
\cite[\S2.9]{6} with $(a,b,c):=(-m,n+1,-m)$. Then $u_1, u_2, u_3$
above ($u_2$ up to a constant factor) are equal to (1), (17), (21),
respectively, in \cite[\S2.9]{6}. However, this is for $u_1$ and $u_3$
only a formal proof, because there occurs a lower parameter in the
hypergeometric function which is a nonpositive integer.
For a rigorous argument for $u_1$,
consider the solution \cite[2.9(1)]{6} of \eqref{26}
first for $(a,b,c):=(-m,n+1,c)$ and then let $c\to-m$. Also, for $u_3$,
consider the solution \cite[2.9(21)]{6} of \eqref{26} first for $c,a:=-b-m-n-1$
with $b$ general, and then let $b\to n+1$.

For the general theory of solving the hypergeometric differential equation
\eqref{26} see \cite[\S2.3]{3}. In general, for fixed $a,b,c$, and on
a simply connected domain in $\CC$ which avoids the singular points
0, 1 (and $\iy$), one can choose two linearly independent solutions
and have the general solution as an arbitrary linear combination of these two
solutions. For a particular solution the coefficients in the linear
combination can be found from (possibly asymptotic) values of the solution
at two of the three singular points. In our case of solutions $u_1$, $u_2$,
$u_3$, given by \eqref{46}--\eqref{47}, the solutions are rational, so they
exist as one-valued functions on $\CC$ (possibly with a pole in 1).
If we would a priori know only that $u_3(x)=Au_1(x)+Bu_2(x)$ then we can
compute $A=1$ by putting $x=0$ and we can compute
$B=1$ by multiplying both sides of the equality by $(1-x)^{n+1}$,
next putting $x=1$, and then using the Chu-Vandermonde identity
\cite[Corollary 2.2.3]{3} for the evaluation of ${}_2F_1(m+1,-n;m+2;1)$.

The case discussed here is the case $m=0$ in Vid\=unas \cite[\S8]{8}
(trivial monodromy group). In this way he obtained \eqref{1} as
the case $m=0$ of the identity at the end of his section 8.

For $\al>0$ and $n,m$ nonnegative integers we will now prove the following
more general
three-term identity:
\begin{multline}
(1-z)^{-n-1} (1-z^{-1})^{-\al}\,
\hyp21{m+1,-\al}{n+m+2}{z^{-1}}\\
=\frac{(n+1)_{m+1}}{(n+\al+1)_{m+1}}\,z^{m+1} (1-z)^{-n-1}\,(1-z^{-1})^{-\al}\,
\hyp21{-n,m+1}{-n-\al}{1-z}\\
+\frac{(m+1)_{n+1}}{(m+\al+1)_{n+1}}\,
\hyp21{-m,n+1}{-m-\al}z\quad\bigl(z\in\CC\backslash[0,1]\,\bigr).
\label{36}
\end{multline}
This formula makes good sense on the indicated domain, since
${}_2F_1(a,b;c;z)$, originally defined as a power series for $|z|<1$,
has a unique analytic continuation to $\CC\backslash[1,\iy)$. Hence,
the ${}_2F_1$ on the left is uniquely defined for $z\notin[0,1]$.
Also, $(1-z^{-1})^\al$ is uniquely defined for $z\notin[0,1]$. The
two ${}_2F_1$'s on the right,
being polynomials, are defined for all $z\in\CC$.
If one wishes, one may rewrite $(1-z)^{-n-1} (1-z^{-1})^{-\al}$ as
$(-1)^{n+1}z^{-n-1}(1-z^{-1})^{-\al-n-1}$.

For $\al\to0$ the identity \eqref{36} tends to the identity \eqref{1}.
In fact, \eqref{36}, which is of the form $u_3=u_1+u_2$ (see \eqref{37}), has
the terms $u_1$, $u_2$ and $u_3$ as solutions
of the differential
equation
\begin{equation}
z(1-z)u''(z)-\bigl((n+2)z+m(1-z)+\al\bigr)u'(z)+m(n+1)u(z)=0,
\label{38}
\end{equation}
i.e., the hypergeometric differential equation \eqref{26} for
$(a,b,c)=(-m,n+1,-m-\al)$.
Indeed, see the solutions \cite[2.9 (18),(1),(14)]{6} of \eqref{38}.
This gives $u_1$ (after substituting \eqref{40}), $u_2$ and $u_3$,
respectively.
For $\al\downarrow 0$, \eqref{38} tends to \eqref{39}, and the solutions
$u_1,u_2,u_3$ of \eqref{38} tend to the solutions $u_1,u_2,u_3$ of
\eqref{39}.

The case of \eqref{36} that $\al$ is a nonnegative integer is essentially
the general case of the identity at the end of section 8 in Vid\=unas \cite{8}.
Just transform the ${}_2F_1$ on the left of \eqref{36} by first reversing
the order of summation in the hypergeometric series and next applying
Pfaff's transformation formula \eqref{48}.

For the proof of \eqref{36}, start with the three-term identity
\begin{multline}
\hyp21{m+1,-\al}{n+m+2}z
=\frac{\Ga(n+m+2)\Ga(n+\al+1)}{\Ga(n+1)\Ga(n+m+\al+2)}\,
z^{-m-1}\,\hyp21{m+1,-n}{-n-\al}{1-z^{-1}}\\
+\frac{\Ga(n+m+2)\Ga(-n-\al-1)}{\Ga(m+1)\Ga(-\al)}\,
z^{-n-1}(1-z)^{n+1+\al}\,\hyp21{n+1,-m}{n+\al+2}{1-z^{-1}}\\
\bigl(z\notin\{0\}\union[1,\iy)\,\bigr),
\label{41}
\end{multline}
see \cite[2.10(4)]{6} (or \cite[(2.3.11)]{3} combined with \eqref{48}).
Note that we don't have to exclude $z\in(-\iy,0)$ in \eqref{41} because
the two ${}_2F_1$'s on the right are terminating.
By \eqref{40} the last hypergeometric function on the right can be
replaced by
\[
\frac{(\al+1)_m}{(n+\al+2)_m}\,\hyp21{n+1,-m}{-m-\al}{z^{-1}}.
\]
In the identity which thus results from \eqref{41}, first replace $z$ by
$z^{-1}$
and next multiply both sides by $(1-z)^{-n-1}(1-z^{-1})^{-\al}$. This yields
\eqref{36}.
\section{A multivariable generalization}
\label{108}
A multivariable generalization \eqref{111} of \eqref{2} was proved in
\cite{16} by generating functions (see the one-variable case of this
proof in section \ref{42}), while a probabilistic proof, immediately
generalizing the one-variable case discussed in Remark \ref{118}, was
indicated in \cite{16} and \cite{13}. We will give in section \ref{110}
a different proof of \eqref{111}, which will generalize the proof of
\eqref{2} in section \ref{17} using the beta integral.

Let us reformulate \eqref{111} in other notation, and let us also give
this identity in homogeneous form.
Let $s:=x_1+\cdots+x_n$. Define
\begin{multline}
f_{a_1,\ldots,a_n}(x_1,\ldots,x_n):=
x_n^{a_n+1}\sum_{k_1=0}^{a_1}\ldots\sum_{k_{n-1}=0}^{a_{n-1}}
\frac{(a_n+1)_{k_1+\cdots+k_{n-1}}}{k_1!\ldots k_{n-1}!}\\
\times
x_1^{k_1}\ldots x_{n-1}^{k_{n-1}}\,s^{a_1+\cdots+a_{n-1}-(k_1+\cdots+k_{n-1})}.
\label{49}
\end{multline}
Then
\begin{equation}
(x_1+\cdots+x_n)^{a_1+\cdots+a_n+1}=\sum_\si
f_{a_{\si(1)},\ldots,a_{\si(n)}}(x_{\si(1)},\ldots,x_{\si(n)}),
\label{55}
\end{equation}
where summation is over all cyclic permutations $\si$
of $1,2,\ldots,n$.
For $x_1+\cdots+x_n=1$ identity \eqref{55} simplifies to
\begin{equation}
1=\sum_\si
f_{a_{\si(1)},\ldots,a_{\si(n)}}(x_{\si(1)},\ldots,x_{\si(n)})
\label{51}
\end{equation}
with
\begin{equation}
f_{a_1,\ldots,a_n}(x_1,\ldots,x_n)=
x_n^{a_n+1}\sum_{k_1=0}^{a_1}\ldots\sum_{k_{n-1}=0}^{a_{n-1}}
\frac{(a_n+1)_{k_1+\cdots+k_{n-1}}}{k_1!\ldots k_{n-1}!}\,
x_1^{k_1}\ldots x_{n-1}^{k_{n-1}}.
\label{52}
\end{equation}
Identity \eqref{51} is a reformulation of \eqref{111} and \eqref{55}
is the homogeneous form of \eqref{51}.
\begin{remark}
For $n=2$ we can write \eqref{51}, \eqref{52} as
\[
1=f_{m,n}(x,1-x)+f_{n,m}(1-x,x)\quad{\rm with}\quad
f_{m,n}(x,1-x)=(1-x)^{n+1}\sum_{k=0}^m\frac{(n+1)_k}{k!}\,x^k.
\]
So we have \eqref{2}. The case $n=3$ of \eqref{51}, \eqref{52} is
also noteworthy as a two-variable analogue of \eqref{2}:
\begin{align}
&1=f_{a,b,c}(x,y,1-x-y)+f_{b,c,a}(y,1-x-y,x)
+f_{c,a,b}(1-x-y,x,y),
\label{53}\\
&f_{a,b,c}(x,y,1-x-y)=(1-x-y)^{c+1}\,
\sum_{k=0}^a\sum_{l=0}^b\frac{(c+1)_{k+l}}{k!\,l!}\,x^ky^l.
\label{54}
\end{align}
Identity \eqref{53} is also given in \cite{16}.

If we divide both sides of \eqref{53} by $(1-x-y)^{c+1}$ then the
resulting identity has the form
\[
(1-x-y)^{-c-1}=f(x,y)+x^{a+1}g(x,y)+y^{b+1}h(x,y),
\]
where $f(x,y)$ is the polynomial consisting of all
terms of the power series of $(1-x-y)^{-c-1}$ which
have degree $\le a$ in $x$ and degree $\le b$ in $y$, while $g(x,y)$ and
$h(x,y)$ are power series in $x$ and $y$. There does not seem to be an
a priori symmetry argument which settles \eqref{53} from this observation.
For instance, if we try to imitate the proof for $n=2$ in
section \ref{33} then we have to consider the homogeneous version
of \eqref{53} given by \eqref{55} for $n=3$:
\[
(x+y+z)^{a+b+c+1}=f_{a,b,c}(x,y,z)+f_{b,c,a}(y,z,x)+f_{c,a,b}(z,x,y).
\]
The problem is where to put on the right the terms
$\ga_{k,l,m}x^ky^lz^m$ ($k+l+m=a+b+c+1$) in the expansion of
$(x+y+z)^{a+b+c+1}$.
Certainly, we can uniquely put all terms $\ga_{k,l,m}x^ky^lz^m$
with $k\le a$, $l\le b$ in $f_{a,b,c}(x,y,z)$, all terms
with $l\le b$, $m\le c$ in $f_{b,c,a}(y,z,x)$,
and all terms
with $m\le c$, $k\le a$ in $f_{c,a,b}(z,x,y)$, but there is no clear rule
where to put a term in which only one of the three inequalities
$k\le a$, $l\le b$, $m\le c$ holds.
\end{remark}
\begin{remark}
Note that we can formally express the \RHS\ of \eqref{54} in terms of
{\em Appell's
hypergeometric function} $F_2$ (see \cite[\S5.7.1]{6}):
\[
f_{a,b,c}(x,y,1-x-y)=(1-x-y)^{c+1}\,
F_2(c+1,-a,-b,-a,-b,x,y).
\]
However, due to the nonpositive integer bottom parameters
we cannot transform this $F_2$ function
by \cite[5.11(8)]{6} similarly
as we transformed
the ${}_2F_1$ function in Remark \ref{56}.

Similarly to the case $n=3$, the multisum on the \RHS\ of \eqref{52}
can be formally written as a
{\em Lauricella hypergeometric function} $F_A$ (see
\cite[Ch.~VII]{10} or
\cite[(8.6.1)]{9}):
\begin{multline}
f_{a_1,\ldots,a_n}(x_1,\ldots,x_n)=
x_n^{a_n+1}
F_A^{(n-1)}(a_n+1,-a_1,\ldots,-a_{n-1},-a_1,\ldots,-a_{n-1};
x_1,\ldots,x_{n-1})\\
(x_1+\cdots+x_n=1).
\label{60}
\end{multline}
Here the $F_A^{(n-1)}$ has to be interpreted as
\begin{equation}
\lim_{(b_1,\ldots,b_{n-1})\to(-a_1,\ldots,-a_{n-1})}
F_A^{(n-1)}(a_n+1,-a_1,\ldots,-a_{n-1},b_1,\ldots,b_{n-1};
x_1,\ldots,x_{n-1}).
\label{61}
\end{equation}
\end{remark}
\section{A PDE associated with the multivariable identity}
\label{109}
{}From the expression in \eqref{60}
of $f_{a_1,\ldots,a_n}(x_1,\ldots,x_{n-1},1-x_1-\cdots-x_{n-1})$
in terms of a Lauricella hypergeometric function, we will derive a
PDE for $f_{a_1,\ldots,a_n}$. Consider first the system of PDE's for the
$F_A^{(n-1)}$ in \eqref{61}, as given in \cite[Ch.~VII, \S XXXIX]{10} 
(for $n=3$ we have Appell's hypergeometric function $F_2$ and then
the system of PDE's is also given in
\cite[5.9(10)]{6}).
After taking the limit for $(b_1,\ldots,b_{n-1})\to(-a_1,\ldots,-a_{n-1})$
we obtain that
\begin{equation}
u(x_1,\ldots,x_{n-1}):=
F_A^{(n-1)}(a_n+1,-a_1,\ldots,-a_{n-1},-a_1,\ldots,-a_{n-1};
x_1,\ldots,x_{n-1})
\label{64}
\end{equation}
satisfies the system of PDE's
\begin{multline*}
x_i(1-x_i)\pa_i^2u-x_i\sum_{j\ne i}x_j\pa_j\pa_i u
-(a_i+(a_n+2)x_i)\pa_i u+a_i\sum_j x_j \pa_j u+(a_n+1)a_iu=0\\
(i=1,\ldots,n-1).
\end{multline*}
Here $\pa_i$ denotes $\pa/\pa x_i$.
The sum of the $n-1$ PDE's equals
\begin{multline}
\sum_i x_i(1-x_i)\pa_i^2u-\sum_{i\ne j}x_ix_j\pa_i\pa_j u
+(a_1+\cdots+a_{n-1}-a_n-2)\sum_i x_i\pa_i u
-\sum_i a_i\pa_i u\\
+(a_n+1)(a_1+\cdots+a_{n-1})u=0.
\label{63}
\end{multline}
\begin{proposition}
The function \eqref{64} is the unique solution, up to a constant factor,
of \eqref{63} which has the form
\begin{equation}
u(x_1,\ldots,x_{n-1})=
\sum_{k_1=0}^{a_1}\ldots\sum_{k_{n-1}=0}^{a_{n-1}}
\ga_{k_1,\ldots,k_{n-1}}x_1^{k_1}\ldots x_{n-1}^{k_{n-1}}.
\label{65}
\end{equation}
\end{proposition}
\Proof
Computation of the \LHS\ of \eqref{63} with
$u:=x_1^{k_1}\ldots x_{n-1}^{k_{n-1}}$ yields
\begin{multline*}
\bigl((a_1-k_1)+\cdots+(a_{n-1}-k_{n-1})\bigr)
(a_n+1+k_1+\cdots+k_{n-1})x_1^{k_1}\ldots x_{n-1}^{k_{n-1}}\\
+\left(\sum_{i=1}^{n-1}k_i(k_i-a_i-1)x_i^{-1}\right)
x_1^{k_1}\ldots x_{n-1}^{k_{n-1}}.
\end{multline*}
It follows that $u$ of the form \eqref{65} satisfies \eqref{63} iff
\begin{multline}
\bigl((a_1-k_1)+\cdots+(a_{n-1}-k_{n-1})\bigr)
(a_n+1+k_1+\cdots+k_{n-1})\ga_{k_1,\ldots,k_{n-1}}\\
+\sum_{i=1}^{n-1}(k_i+1)(k_i-a_i)\ga_{k_1,\ldots,k_i+1,\ldots,k_{n-1}}=0.
\label{66}
\end{multline}
Give some value to $\ga_{a_1,\ldots,a_{n-1}}$.
Then we see from \eqref{66}
by downward induction with respect to $k_1+\cdots+k_{n-1}$
that all coefficients $\ga_{k_1,\ldots,k_{n-1}}$ with $0\le k_i\le a_i$
($i=1,\ldots,n-1$) are uniquely determined by this initial value.

In passing we see that \eqref{66} is satisfied by
\[
\ga_{k_1,\ldots,k_{n-1}}:=
\frac{(a_n+1)_{k_1+\cdots+k_{n-1}}}{k_1!\,\ldots k_{n-1}!}\,.
\]
Thus we have also proved from scratch that $u$ given by \eqref{64} satisfies
\eqref{63}.
\qed
\bPP
By some computation, we see that
\begin{align*}
v(x_1,\ldots,x_{n-1})&:=
(1-x_1-\cdots-x_{n-1})^{a_n+1}u(x_1,\ldots,x_{n-1})\\
&\;=f_{a_1,\ldots,a_n}(x_1,\ldots,x_{n-1},1-x_1-\cdots-x_{n-1})
\end{align*}
satisfies the PDE
\begin{equation}
\sum_{i=1}^{n-1}x_i(1-x_i)\pa_i^2v-2\sum_{i<j}x_ix_j\pa_i\pa_jv+
\sum_{i=1}^{n-1}\Bigl((a_1+\cdots+a_n)x_i-a_i\Bigr)\pa_iv=0.
\label{62}
\end{equation}
Clearly, the function $v:=1$ satisfies \eqref{62}.
Furthermore, by straightforward computations we see: If $v$
satisfies \eqref{62} then the function
$(x_1,\ldots,x_{n-1})\allowbreak
\mapsto v(1-x_1-\cdots-x_{n-1},x_2,\ldots,x_{n-1})$
satisfies \eqref{62} with $a_1$ and $a_n$ interchanged.
Thus we have proved:
\begin{theorem}
For all permutations $\si$ of $1,2,\ldots,n$ the functions
\[
(x_1,\ldots,x_{n-1})\mapsto
f_{a_{\si(1)},\ldots,a_{\si(n)}}(x_{\si(1)},\ldots,x_{\si(n)})\qquad
(x_1+\cdots+x_n=1)
\]
are solutions of \eqref{62}. Up to a constant factor they are the unique
solutions of \eqref{62} of the form
\begin{equation*}
v(x_1,\ldots,x_{n-1})=x_{\si(n)}^{a_{\si(n)}+1}\,
\sum_{k_1=0}^{a_{\si(1)}}\ldots\sum_{k_{n-1}=0}^{a_{\si(n-1)}}
\ga_{k_1,\ldots,k_{n-1}\phantom{k^2}\!\!\!\!\!}
x_{\si(1)}^{k_1}\ldots x_{\si(n-1)}^{k_{n-1}}\qquad
(x_1+\cdots+x_n=1).\qquad\qed
\end{equation*}
\label{67}
\end{theorem}

Now consider a solution $v(x_1,\ldots,x_{n-1})$ of \eqref{62}, let
$x_n$ be a variable independent of $x_1,\ldots,x_{n-1}$, and let $\phi$
be an arbitrary function of that new variable. Then trivially \eqref{62}
holds with $v$ replaced by $\phi(x_n)v(x_1,\ldots,x_{n-1})$.
Now pass in this PDE to new variables $y_1,\ldots,y_n$ by
\[
x_1=\frac{y_1}{y_1+\cdots+y_n}\,,\;\ldots,\;
x_{n-1}=\frac{y_{n-1}}{y_1+\cdots+y_n}\,,\quad
x_n=y_1+\cdots+y_n,
\]
or equivalently,
\[
y_1=x_1x_n,\;\ldots,\;y_{n-1}=x_{n-1}x_n,\quad
y_n=(1-x_1-\cdots-x_{n-1})x_n.
\]
Then we obtain that the function
$w(y_1,\ldots,y_n):=\phi(y_1+\cdots+y_n)\,
v\bigl(\frac{y_1}{y_1+\cdots+y_n},\ldots,\frac{y_{n-1}}{y_1+\cdots+y_n}\bigr)$
satisfies the PDE
\begin{equation}
\sum_{i=1}^n y_i(y_1+\cdots+y_n-y_i)\pa_i^2w-2\sum_{i<j}y_iy_j\pa_i\pa_jw
+\sum_{i=1}^n\Bigl((a_1+\cdots+a_n)y_i-a_i(y_1+\cdots+y_n)\Bigr)\pa_iw=0,
\label{68}
\end{equation}
where $\pa_i$ denotes $\pa/\pa y_i$.
Thus by \eqref{49}, \eqref{52} and
Theorem \ref{67} we have proved in particular:
\begin{theorem}
The function
$w(y_1,\ldots,y_n):=f_{a_1,\ldots,a_n}(y_1,\ldots,y_n)$, defined by
\eqref{49}, satisfies \eqref{68}.
Similarly, the functions
$w(y_1,\ldots,y_n):=
f_{a_{\si(1)},\ldots,a_{\si(n)}}(y_{\si(1)},\ldots,y_{\si(n)})$
satisfy \eqref{68},
where $\si$ is a permutation of $1,\ldots,n$.
Also, the function $w(y_1,\ldots,y_n):=(y_1+\cdots+y_n)^{a_1+\cdots+a_n+1}$
satisfies \eqref{68}.
\end{theorem}
\section{Splitting up Dirichlet's multivariable beta integral}
\label{110}
Just as \eqref{1} can be obtained by splitting a beta integral into two
parts and evaluating the resulting incomplete beta integrals, we can
prove and interprete \eqref{51} by splitting Dirichlet's $(n-1)$-dimensional
beta integral with nonnegative integer exponents
into $n$ parts. For convenience, we will work here with
an $n$-dimensional beta integral.

Let $\De_n$ be the simplex in $\RR^n$ wich has vertices $0$ and the standard
basis vectors $e_1,\ldots,e_n$.
Let $a_1,\ldots,a_n,b$ be complex numbers with real part $>-1$. Then
Dirichlet's integral is as follows.
\begin{equation}
I_{a_1,\ldots,a_{n+1}}:=
\int_{\De_n}t_1^{a_1}\ldots t_n^{a_n}\,(1-t_1-\cdots-t_n)^{a_{n+1}}\,
dt_1\ldots dt_n
=\frac{\Ga(a_1+1)\ldots\Ga(a_{n+1}+1)}
{\Ga(a_1\cdots+a_{n+1}+n+1)}\,,
\label{101}
\end{equation}
see \cite[Theorem 1.8.6]{3} or \cite[Exercise 7.2.6]{11} for a straightforward
proof, and \cite{12} for its history.
Note that $I_{a_1,\ldots,a_{n+1}}$ is symmetric in $a_1,\ldots,a_{n+1}$.

Now take $x=(x_1,\ldots,x_n)$ within $\De_n$ and let
$\De_n^{(i)}(x)$ ($i=1,\ldots,n+1$)
denote the simplex in $\RR^n$ which has a vertex $x$ and $n$
further vertices taken from
$0,e_1,\ldots,e_n$ where $e_i$ is deleted if $i=1,\ldots,n$
and 0 is deleted if $i=n+1$.
Require that $a_1,\ldots,a_{n+1}$ are nonnegative integers. Define
\begin{equation}
I_{a_1,\ldots,a_{n+1}}^{(i)}(x):=
\int_{\De_n^{(i)}(x)}t_1^{a_1}\ldots t_n^{a_n}\,(1-t_1-\cdots-t_n)^{a_{n+1}}\,
dt_1\ldots dt_n.
\label{102}
\end{equation}
For any $(y_1,\ldots,y_n)\in\RR^n$ put $y_{n+1}:=1-y_1-\cdots-y_n$.
Then, for any permutation $\si$ of $1,2,\ldots,n+1$ (i.e., $\si\in S_{n+1}$), 
the map
$(y_1,\ldots,y_n)\mapsto(y_{\si(1)},\ldots,y_{\si(n)})$ is a diffeomorphism
of $\De_n$ with Jacobian having absolute value 1. Thus we have
the identity
\begin{equation}
I_{a_1,\ldots,a_{n+1}}^{(i)}(x_1,\ldots,x_n)=
I_{a_{\si(1)},\ldots,a_{\si(n+1)}}^{(\si^{-1}(i))}
(x_{\si(1)},\ldots,x_{\si(n)})\qquad(\si\in S_{n+1}).
\label{105}
\end{equation}
Define
\begin{equation}
f_{a_1,\ldots,a_{n+1}}(x):=(1-x_1-\cdots-x_n)^{a_{n+1}+1}\,
\sum_{k_1=0}^{a_1}\ldots\sum_{k_n=0}^{a_n}
\frac{(a_{n+1}+1)_{k_1+\cdots+k_n}}
{k_1!\,\ldots k_n!}\,x_1^{k_1}\ldots x_n^{k_n},
\label{103}
\end{equation}
i.e., \eqref{49} with $n$ replaced by $n+1$ and with
$x_{n+1}:=1-x_1-\cdots-x_n$ omitted in the argument.
We have the symmetry
\begin{equation}
f_{a_1,\ldots,a_{n+1}}(x_1,\ldots,x_n)
=f_{a_{\si(1)},\ldots,a_{\si(n)},a_{n+1}}(x_{\si(1)},\ldots,x_{\si(n)})\qquad
(\si\in S_n).
\end{equation}
\begin{proposition}
For nonnegative integers $a_1,\ldots,a_{n+1}$ and for $x$ within $\De_n$
we have:
\begin{equation}
\frac{I_{a_1,\ldots,a_{n+1}}^{(n)}(x)}{I_{a_1,\ldots,a_{n+1}}}=
f_{a_{n+	1},a_1,\ldots,a_n}(1-x_1-\cdots-x_n,x_1,\ldots,x_{n-1}).
\label{104}
\end{equation}
\end{proposition}
\Proof
For convenience put $b:=a_{n+1}$ and $x':=(x_1,\ldots,x_{n-1})$.
Then
\begin{align*}
&I_{a_1,\ldots,a_n,b}^{(n)}(x)=\\
&\int_0^{x_n}t_n^{a_n}\left(\int_{(t_n/x_n)x'+(1-(t_n/x_n))\De_{n-1}}
t_1^{a_1}\ldots t_{n-1}^{a_{n-1}}\,(1-t_1-\cdots-t_n)^b\,
dt_1\ldots dt_{n-1}\right)dt_n\\
\noalign{\allowbreak}
&=x_n^{a_n+1}\int_0^1 s^{a_n}\left(
\int_{sx'+(1-s)\De_{n-1}}
t_1^{a_1}\ldots t_{n-1}^{a_{n-1}}\,(1-t_1-\cdots-t_{n-1}-sx_n)^b\,
dt_1\ldots dt_{n-1}\right)ds\\
\noalign{\allowbreak}
&=x_n^{a_n+1}\int_0^1 s^{a_n}\,(1-s)^{n-1}\Biggl(
\int_{\De_{n-1}}(sx_1+(1-s)s_1)^{a_1}\ldots(sx_{n-1}+(1-s)s_{n-1})^{a_{n-1}}\\
&\qquad\qquad\qquad
\times\bigl(s(1-x_1-\cdots-x_n)+(1-s)(1-s_1-\cdots-s_{n-1})\bigr)^b\,
ds_1\ldots ds_{n-1}\Biggr)ds\\
\noalign{\allowbreak}
&=x_n^{a_n+1}\sum_{k_1=0}^{a_1}\ldots\sum_{k_{n-1}=0}^{a_{n-1}}\;\sum_{l=0}^b
\binom{a_1}{k_1}\ldots\binom{a_{n-1}}{k_{n-1}}\binom bl
x_1^{k_1}\ldots x_{n-1}^{k_{n-1}}(1-x_1-\cdots-x_n)^l\\
&\qquad\qquad\qquad\times\int_0^1 s^{a_n+k_1+\cdots+k_{n-1}+l}\,
(1-s)^{a_1+\cdots+a_{n-1}+b+n-1-k_1-\cdots-k_{n-1}-l}\,ds\\
&\qquad\qquad\qquad\times
\int_{\De_{n-1}}s_1^{a_1-k_1}\ldots s_{n-1}^{a_{n-1}-k_{n-1}}\,
(1-s_1-\cdots-s_{n-1})^{b-l}\,ds_1\ldots ds_{n-1}\\
\noalign{\allowbreak}
&=x_n^{a_n+1}\sum_{k_1=0}^{a_1}\ldots\sum_{k_{n-1}=0}^{a_{n-1}}\;\sum_{l=0}^b
\binom{a_1}{k_1}\ldots\binom{a_{n-1}}{k_{n-1}}\binom bl
x_1^{k_1}\ldots x_{n-1}^{k_{n-1}}\,(1-x_1-\cdots-x_n)^l\\
&\qquad\qquad\times
\frac{\Ga(a_n+k_1+\cdots+k_{n-1}+l+1)
\Ga(a_1+\cdots+a_{n-1}+b+n-k_1-\cdots-k_{n-1}-l)}
{\Ga(a_1+\cdots+a_n+b+n+1)}\\
&\qquad\qquad\times
\frac{\Ga(a_1-k_1+1)\ldots\Ga(a_{n-1}-k_{n-1}+1)\Ga(b-l+1)}
{\Ga(a_1+\cdots+a_{n-1}+b+n-k_1-\cdots-k_{n-1}-l)}\\
\noalign{\allowbreak}
&=\frac{a_1!\,\ldots a_n!\,b!}{(a_1+\cdots+a_n+b+n)!}\\
&\qquad\qquad\times x_n^{a_n+1}\,
\sum_{k_1=0}^{a_1}\ldots\sum_{k_{n-1}=0}^{a_{n-1}}\;\sum_{l=0}^b
\frac{(a_n+1)_{k_1+\cdots+k_{n-1}+l}}{k_1!\,\ldots k_{n-1}!\,l!}\,
x_1^{k_1}\ldots x_{n-1}^{k_{n-1}}\,(1-x_1-\cdots-x_n)^l\\
\noalign{\allowbreak}
\quad&=I_{a_1,\ldots,a_n,b}\,
f_{b,a_1,\ldots,a_{n-1},a_n}(1-x_1-\cdots-x_n,x_1,\ldots,x_{n-1}).
\qquad\qquad\qquad\qquad\qquad\qquad\quad\qed
\end{align*}
\begin{theorem}
\label{106}
Let $a_1,\ldots,a_{n+1}$ and $x$ as before. Let $i\in\{1,2,\ldots,n+1\}$
and let $\si$ be the cyclic permutation of $1,\ldots,n+1$ which
sends $n$ to $i$. Then
\begin{equation}
\frac{I_{a_1,\ldots,a_{n+1}}^{(i)}(x)}{I_{a_1,\ldots,a_{n+1}}}
=f_{a_{\si(n+1)},a_{\si(1)},\ldots,a_{\si(n)}}
(x_{\si(n+1)},x_{\si(1)},\ldots,x_{\si(n-1)}).
\end{equation}
\end{theorem}
\Proof
By \eqref{105} we have
$\displaystyle
\frac{I_{a_1,\ldots,a_{n+1}}^{(i)}(x)}{I_{a_1,\ldots,a_{n+1}}}
=\frac{I_{a_{\si(1)},\ldots,a_{\si(n+1)}}^{(n)}(x_{\si(1)},\ldots,x_{\si(n)})}
{I_{a_{\si(1)},\ldots,a_{\si(n+1)}}}$.
Now apply \eqref{104}.
\qed
\bPP
We have the obvious identity
\begin{equation}
1=\sum_{i=1}^{n+1}
\frac{I_{a_1,\ldots,a_{n+1}}^{(i)}(x)}{I_{a_1,\ldots,a_{n+1}}}\quad
(x\in\De_n,\;a_1,\ldots,a_{n+1}\in\CC,\;\Re a_1,\ldots,\Re a_{n+1}>-1).
\label{124}
\end{equation}
By Theorem \ref{106} this
is for nonnegative integers $a_1,\ldots,a_{n+1}$
equivalent with \eqref{51} (with $n$ replaced by $n+1$ and with the
functions $f_{a_1,\ldots,a_{n+1}}$ defined by \eqref{52}). In the
general case of \eqref{124} we get an extension of \eqref{51} for
non-integer $a_1,\ldots,a_{n+1}$, just as we discussed in the one-variable
case in section \ref{19}. The uniqueness of the extension if the terms
satisfy estimates as in Carlson's theorem, as discussed there, also holds here.
\begin{remark}
It is an interesting question (but for us a nontrivial open problem)
to find an elegant looking evaluation of \eqref{102} which is valid for
all complex $a_1,\ldots,a_{n+1}$ with real part $>-1$, and which would
generalize the evaluation \eqref{35} of the incomplete beta function.
This would also give an $n$-variable generalization of \eqref{23},
i.e., of a limit case of Pfaff's identity \eqref{40}.
\end{remark}

\quad\\
\begin{footnotesize}
\begin{quote}
T.~H. Koornwinder, Korteweg-de Vries Institute, University of Amsterdam,\\
Plantage Muidergracht 24, 1018 TV Amsterdam, The Netherlands;\\
email: {\tt T.H.Koornwinder@uva.nl}
\bLP
M.~J. Schlosser,
Fakult\"at f\"ur Mathematik, Universit\"at Wien,\\
Nordbergstrasse 15,
A-1090 Vienna, Austria;\\
email: {\tt michael.schlosser@univie.ac.at}
\end{quote}
\end{footnotesize}

\end{document}